\documentclass[11pt]{article}

\usepackage{amssymb,amsthm,amsmath}
\usepackage[cp1250]{inputenc}
\usepackage{amssymb}
\usepackage[a4paper, left=3cm, right=3.5cm, textwidth=14.5cm]{geometry}

\newcommand{\p}{\mathbb{P}}
\newcommand{\N}{\mathbb{N}}
\newcommand{\E}{\mathbb{E}}
\newcommand{\Var}{\mathrm{Var}}
\newcommand{\R}{\mathbb{R}}

\newtheorem{theorem}{Theorem}
\newtheorem{lemma}{Lemma}

\newtheorem{defi}{Definition}

\newcommand*{\ind}[1]{\mathbf{1}_{\{#1\}}}

\author{Rados{\l}aw Adamczak\thanks{
Research partially supported by MEiN Grant 1 PO3A 012 29.}\\ \\
email: R.Adamczak@impan.gov.pl
}

\title{A tail inequality for suprema of unbounded empirical processes with applications to Markov chains}

\begin{document}
\maketitle

\begin{abstract}
We present a tail inequality for suprema of empirical processes
generated by variables with finite $\psi_\alpha$ norms and apply it to
some geometrically ergodic Markov chains to derive similar
estimates for empirical processes of such chains, generated by
bounded functions. We also obtain a bounded difference inequality
for symmetric statistics of such Markov chains.\\



\vskip0.5cm

\end{abstract}
\newpage
\section{Introduction}

Let us consider a sequence $X_1,X_2,\ldots,X_n$ of random
variables with values in a measurable space
$(\mathcal{S},\mathcal{B})$ and a countable class of measurable
functions $f\colon \mathcal{S} \to \R$. Define moreover the random
variable
\begin{displaymath}
Z = \sup_{f\in \mathcal{F}} |\sum_{i=1}^n f(X_i)|.
\end{displaymath}

In recent years a lot of effort has been devoted to describing the
behaviour, in particular concentration properties of the variable
$Z$ under various assumptions on the sequence $X_1,\ldots,X_n$ and
the class $\mathcal{F}$. Classically, one considers the case of
i.i.d. or independent random variables $X_i$'s and uniformly
bounded classes of functions, although there are also results for
unbounded functions or sequences of variables satisfying some
mixing conditions.

The aim of this paper is to present tail inequalities for the
variable $Z$ under two different types of assumptions, relaxing
the classical conditions.

In the first part of the article we consider the case of
independent variables and unbounded functions (satisfying however
some integrability assumptions). The main result of this part is
Theorem \ref{tail_estimate}, presented in Section
\ref{independent}.

In the second part we keep the assumption of uniform boundedness
of the class $\mathcal{F}$ but relax the condition on the
underlying sequence of variables, by considering a class of Markov
chains, satisfying classical small set conditions with
exponentially integrable regeneration times. If the small set
assumption is satisfied for the one step transition kernel, the
regeneration technique for Markov chains together with the results
for independent variables and unbounded functions allow us to
derive tail inequalities for the variable $Z$ (Theorem
\ref{empirical_Markov}, presented in Section
\ref{subsectionMarkov}).

In a more general situation, when the small set assumption is
satisfied only by the $m$-skeleton chain, our results are
restricted to sums of real variables, i.e. to the case of
$\mathcal{F}$ being a singleton (Theorem \ref{sums_Markov}).

Finally, in Section \ref{boundeddifference}, using similar
arguments, we derive a bounded difference type inequality for
Markov chains, satisfying the same small set assumptions.

We will start by describing known results for bounded classes of
functions and independent random variables, beginning with the
celebrated Talagrand's inequality. They will serve us both as
tools and as a point of reference for presenting our results.

\subsection{Talagrand's concentration inequalities}

In the paper \cite{T}, Talagrand proved the following inequality
for empirical processes.

\begin{theorem}[Talagrand, \cite{T}] \label{Talagrand} Let $X_1,\ldots,X_n$ be independent random
variables with values in a measurable space
$(\mathcal{S},\mathcal{B})$ and let $\mathcal{F}$ be a countable
class of measurable functions $f\colon \mathcal{S} \to \R$, such
that $\|f\|_\infty \le a < \infty$ for every $f \in \mathcal{F}$.
Consider the random variable $Z = \sup_{f \in
\mathcal{F}}\sum_{i=1}^n f(X_i)$. Then for all $t \ge 0$,
\begin{align}\label{Talagrand_1}
\p(Z \ge \E Z + t) \le K\exp\Big(-\frac{1}{K}\frac{t}{a}\log(1 +
\frac{ta}{V})\Big),
\end{align}
where $V = \E\sup_{f\in \mathcal{F}} \sum_{i=1}^n f(X_i)^2$ and
$K$ is an absolute constant. In consequence, for all $t \ge 0$,
\begin{align}\label{Talagrand_2}
\p(Z \ge \E Z + t) \le K_1\exp\Big(-\frac{1}{K_1}\frac{t^2}{V +
at}\Big)
\end{align}
for some universal constant $K_1$. Moreover, the above
inequalities hold, when replacing $Z$ by $-Z$.
\end{theorem}

Inequalities \ref{Talagrand_1} and \ref{Talagrand_2} may be
considered functional versions of respectively Bennett's and
Bernstein's inequalities for sums of independent random variables
and similarly as in the classical case, one of them implies the
other. Let us note, that Bennett's inequality recovers both the
subgaussian and Poisson behaviour of sums of independent random
variables, corresponding to classical limit theorems, whereas
Bernstein's inequality recovers the subgaussian behaviour for
small values and exhibits exponential behaviour for larger values
of $t$.

The above inequalities proved to be a very important tool in
infinite dimensional probability, machine learning and
M-estimation. They drew considerable attention resulting in
several simplified proofs and different versions. In particular,
there has been a series of papers, starting from the work by
Ledoux \cite{L1}, exploring concentration of measure for empirical
processes with the use of logarithmic Sobolev inequalities with
discrete gradients. The first explicit constants were obtained by
Massart \cite{Mass}, who proved in particular the following

\begin{theorem}[Massart, \cite{Mass}] \label{Massart} Let $X_1,\ldots,X_n$ be independent random
variables with values in a measurable space
$(\mathcal{S},\mathcal{B})$ and let $\mathcal{F}$ be a countable
class of measurable functions $f\colon \mathcal{S} \to \R$, such
that $\|f\|_\infty \le a < \infty$ for every $f \in \mathcal{F}$.
Consider the random variable $Z = \sup_{f \in
\mathcal{F}}|\sum_{i=1}^n f(X_i)|$. Assume moreover that for all
$f\in\mathcal{F}$ and all $i$, $\E f(X_i) = 0$ and let $\sigma^2 =
\sup_{f\in\mathcal{F}} \sum_{i=1}^n \E f(X_i)^2$. Then for all
$\eta
> 0$ and $t \ge 0$,
\begin{displaymath}
\p(Z \ge \ (1+\eta)\E Z + \sigma\sqrt{2K_1t} + K_2(\eta)a t) \le
e^{-t}
\end{displaymath}
and
\begin{displaymath}
\p(Z \le \ (1-\eta)\E Z - \sigma\sqrt{2K_3t} - K_4(\eta)a t) \le
e^{-t},
\end{displaymath}
where $K_1 = 4$, $K_2(\eta) = 2.5 + 32/\eta$, $K_3 = 5.4$,
$K_4(\eta) = 2.5 + 43.2/\eta$.
\end{theorem}

Similar, more refined results were obtained subsequently by
Bousquet \cite{Bou} and Klein and Rio \cite{KR}. The latter
article contains an inequality for suprema of empirical processes
with the best known constants.
\begin{theorem}[Klein, Rio, \cite{KR}, Theorems 1.1., 1.2]\label{KleinRio}
Let $X_1,X_2,\ldots,X_n$ be independent random variables with
values in a measurable space $(S,\mathcal{B})$ and let
$\mathcal{F}$ be a countable class of measurable functions $f
\colon \mathcal{S}\to [-a,a]$, such that for all $i$, $\E f(X_i)
=0$. Consider the random variable
\begin{displaymath}
Z =\sup_{f\in \mathcal{F}} \sum_{i=1}^n f(X_i).
\end{displaymath}
Then, for all $t \ge 0$,
\begin{displaymath}
\p(Z \ge \E Z +t) \le \exp\Big(-\frac{t^2}{2(\sigma^2 + 2a\E Z) +
3at}\Big)
\end{displaymath}
and
\begin{displaymath}
\p(Z \le \E Z - t) \le \exp\Big(-\frac{t^2}{2(\sigma^2 + 2a\E Z) +
3at}\Big),
\end{displaymath}
where
\begin{displaymath}
\sigma^2 = \sup_{f\in \mathcal{F}}\sum_{i=1}^n \E f(X_i)^2.
\end{displaymath}
\end{theorem}

The reader may notice, that contrary to the original Talagrand's
result, estimates of Theorem \ref{Massart} and \ref{KleinRio} use
rather the 'weak' variance $\sigma^2$ than the 'strong' parameter
$V$ of Theorem \ref{Talagrand}. This stems from several reasons,
e.g. the statistical relevance of parameter $\sigma$ and analogy
with the concentration of Gaussian processes (which by CLT, in the
case of Donsker classes of functions correspond to the limiting
behaviour of empirical processes). One should also note, that by
the contraction principle we have $\sigma^2 \le V \le \sigma^2 +
16a\E Z$ (see \cite{LT}, Lemma 6.6). Thus, usually one would
like to describe the subgaussian behaviour of the variables $Z$
rather in terms of $\sigma$, however the price to be paid is the
additional summand of the form $\eta\E Z$. Let us also remark, that if one
does not pay attention to constants, inequalities
presented in Theorems \ref{Massart} and \ref{KleinRio} follow from
Talagrand's inequality just by the aforementioned estimate $V \le
\sigma^2 + 16a\E Z$ and the inequality between the geometric and
the arithmetic mean (in the case of Theorem \ref{Massart}).

\subsection{Notation, basic definitions}

 In the
article, by $K$ we will denote universal constants and by
$C(\alpha,\beta), K_\alpha$ -- constants depending only on
$\alpha,\beta$ or only on $\alpha$ resp. (where $\alpha, \beta$
are some parameters). In both cases the values of constants may
change from line to line.

We will also use the classical definition of (exponential) Orlicz
norms.

\begin{defi}\label{Orlicz_norm} For $\alpha > 0$, define the function $\psi_\alpha \colon \R_+ \to \R_+$ with the formula
$\psi_\alpha(x) = \exp(x^\alpha) - 1$. For a random variable $X$,
define also the Orlicz norm
\begin{displaymath}
\|X\|_{\psi_\alpha} = \inf\{\lambda > 0 \colon
\E\psi_\alpha(|X|/\lambda) \le 1\}.
\end{displaymath}
\end{defi}

Let us also note a basic fact that we will use in the sequel,
namely that by Chebyshev's inequality, for $t\ge 0$,
\begin{displaymath}
\p(|X| \ge t) \le
2\exp\Big(-\Big(\frac{t}{\|X\|_{\psi_\alpha}}\Big)^{\alpha}\Big).
\end{displaymath}

\paragraph{Remark}

For $\alpha < 1$ the above definition does not give a norm but
only a quasi-norm. It can be fixed by changing the function
$\psi_\alpha$ near zero, to make it convex (which would give an
equivalent norm). It is however widely accepted in literature to
use the word norm also for the quasi-norm given by our definition.

\section{Tail inequality for suprema of empirical processes corresponding to classes of unbounded functions \label{independent}}

\subsection{The main result for the independent case}
We will now formulate our main result in the setting of
independent variables, namely tail estimates for suprema of
empirical processes under the assumption that the summands have
finite $\psi_\alpha$ Orlicz norm.

\begin{theorem}\label{tail_estimate}
Let $X_1,\ldots,X_n$ be independent random variables with values
in a measurable space $(\mathcal{S},\mathcal{B})$ and let
$\mathcal{F}$ be a countable class of measurable functions
$f\colon \mathcal{S} \to \R$. Assume that for every $f\in
\mathcal{F}$ and every $i$, $\E f(X_i) = 0$ and for some $\alpha
\in (0,1]$ and all $i$, $\|\sup_{f}| f(X_i)|\|_{\psi_\alpha} <
\infty$. Let
\begin{displaymath}\label{tail_estimate_eq}
Z = \sup_{f\in\mathcal{F}} |\sum_{i=1}^n f(X_i)|.
\end{displaymath}
Define moreover
\begin{displaymath}
\sigma^2 = \sup_{f\in\mathcal{F}} \sum_{i=1}^n\E f(X_i)^2.
\end{displaymath}

Then, for all $0< \eta < 1$ and $\delta
> 0$, there exists a constant $C =C(\alpha,\eta,\delta)$, such
that for all $t \ge 0$,
\begin{align*}
\p(Z &\ge (1+\eta)\E Z + t) \\
&\le \exp\Big(-\frac{t^2}{2(1+\delta)\sigma^2}\Big) +
3\exp\Big(-\Big(\frac{t}{C\|\max_{i}\sup_{f\in\mathcal{F}}|f(X_i)|\|_{\psi_\alpha}}\Big)^{\alpha}\Big)
\end{align*}
and
\begin{align*}
\p(Z &\le (1-\eta)\E Z - t) \\
&\le \exp\Big(-\frac{t^2}{2(1+\delta)\sigma^2}\Big) +
3\exp\Big(-\Big(\frac{t}{C\|\max_{i}\sup_{f\in\mathcal{F}}|f(X_i)|\|_{\psi_\alpha}}\Big)^{\alpha}\Big).
\end{align*}
\end{theorem}

\paragraph{Remark}
The above theorem may be thus considered a counterpart of
Massart's result (Theorem \ref{Massart}). It is written in a
slightly different manner, reflecting the use of Theorem
\ref{KleinRio} in the proof, but it is easy to see that if one
disregards the constants, it yields another version in flavour of
inequalities presented in Theorem \ref{Massart}.

Let us note that some weaker (e.g. not recovering the proper power
$\alpha$ in the subexponential decay of the tail) inequalities may
be obtained by combining the Pisier inequality (see (\ref{Pisier})
below) with moment estimates for empirical processes proven by
Gin\'{e}, Lata{\l}a and Zinn \cite{GLZ} and later obtained by a
different method also by Bousquet, Boucheron, Lugosi and Massart
\cite{BBLM}. These moment estimates first appeared in the context
of tail inequalities for $U$-statistics and were later used in
statistics, in model selection. They are however also of
independent interest as extensions of classical Rosenthal's
inequalities for $p$-th moments of sums of independent random
variables (with the dependence on $p$ stated explicitly).

\paragraph{} The proof of Theorem \ref{tail_estimate} is a compilation of the classical
Hoffman-J{\o}rgensen inequality with Theorem \ref{KleinRio} and
another deep result due to Talagrand.

\begin{theorem}[Ledoux, Talagrand, \cite{LT}, Theorem 6.21. p. 172]\label{Tal_psi} In the setting of Theorem
\ref{tail_estimate}, we have
\begin{align*}
\|Z\|_{\psi_\alpha} \le K_\alpha\Big(\|Z\|_1 + \Big\|\max_i\sup_f
|f(X_i)|\Big\|_{\psi_\alpha}\Big).
\end{align*}
\end{theorem}

We will also need the following corollary to Theorem
\ref{KleinRio}, which was derived in \cite{EL}. Since the proof is
very short we will present it here for the sake of completeness

\begin{lemma}\label{CLT_type}
In the setting of Theorem \ref{KleinRio}, for all $0 < \eta \le
1$, $\delta > 0$ there exists a constant $C = C(\eta,\delta)$,
such that for all $t \ge 0$,
\begin{align*}
\p(Z \ge (1+\eta)\E Z + t) \le
\exp\Big(-\frac{t^2}{2(1+\delta)\sigma^2}\Big) +
\exp\Big(-\frac{t}{Ca}\Big)
\end{align*}
and
\begin{align*}
\p(Z \le (1-\eta)\E Z - t) \le
\exp\Big(-\frac{t^2}{2(1+\delta)\sigma^2}\Big) +
\exp\Big(-\frac{t}{Ca}\Big).
\end{align*}
\end{lemma}

\begin{proof}
It is enough to notice that for all $\delta > 0$,
\begin{align*}
\exp\Big(-\frac{t^2}{2(\sigma^2 + 2a\E Z) + 3at}\Big) \le&
\exp\Big(-\frac{t^2}{2(1+ \delta)\sigma^2}\Big) \\
&+ \exp\Big(-\frac{t^2}{(1 + \delta^{-1})(4a\E Z + 3ta)}\Big)
\end{align*}
and use this inequality together with Theorem \ref{KleinRio} for $t + \eta\E Z$ instead of $t$, which
gives $C = (1 + 1/\delta)(3 + 2\eta^{-1})$.
\end{proof}

\begin{proof}[Proof of Theorem \ref{tail_estimate}]

Without loss of generality we may and will assume that

\begin{align}\label{thewlogrestriction}
t/\|\max_{1\le i\le n}\sup_{f\in\mathcal{F}}
|f(X_i)|\|_{\psi_\alpha} > K(\alpha,\eta,\delta),
\end{align}
otherwise we can make the theorem trivial by choosing the
constant $C = C(\eta,\delta,\alpha)$ to be large enough. The
conditions on the constant $K(\alpha,\eta,\delta)$ will be imposed
later on in the proof.

Let $\varepsilon = \varepsilon(\delta) > 0$ (its value will be
determined later) and for all $f \in \mathcal{F}$ consider the
truncated functions $f_1(x) = f(x)\ind{\sup_{f\in
\mathcal{F}}|f(x)| \le \rho}$ (the truncation level $\rho$ will
also be fixed later). Define also functions $f_2(x) = f(x) -
f_1(x) = f(x)\ind{\sup_{f\in\mathcal{F}}|f(x)| > \rho}$. Let
$\mathcal{F}_i = \{f_i\colon f\in \mathcal{F}\}$, $i = 1,2$.

We have
\begin{align}\label{firsttriangle}
Z = \sup_{f\in \mathcal{F}}|\sum_{i=1}^n f(X_i)| \le&
\sup_{f_1\in\mathcal{F}_1}|\sum_{i=1}^n (f_1(X_i) - \E f_1(X_i))|\nonumber
\\
&+ \sup_{f_2\in\mathcal{F}_2}|\sum_{i=1}^n (f_2(X_i) - \E
f_2(X_i))|\
\end{align} and
\begin{align}\label{secondtriangle}
Z \ge& \sup_{f_1\in\mathcal{F}_1}|\sum_{i=1}^n (f_1(X_i) - \E
f_1(X_i))|\nonumber
\\
&- \sup_{f_2\in\mathcal{F}_2}|\sum_{i=1}^n (f_2(X_i) - \E
f_2(X_i))|
\end{align}
where we used the fact that $\E f_1(X_i) + \E f_2(X_i) = 0$ for
all $f \in \mathcal{F}$.

Similarly, by Jensen's inequality, we get
\begin{align}\label{etriangle}
\E \sup_{f_1\in\mathcal{F}_1} |\sum_{i=1}^n (f_1(X_i)- \E
f_1(X_i))| - 2\E \sup_{f_2\in\mathcal{F}_2} |\sum_{i=1}^n
f_2(X_i)| \le \E Z \nonumber\\
\le \E \sup_{f_1\in\mathcal{F}_1} |\sum_{i=1}^n (f_1(X_i)- \E
f_1(X_i))| + 2\E\sup_{f_2\in\mathcal{F}_2}|\sum_{i=1}^n f_2(X_i)|.
\end{align}
Denoting \begin{displaymath} A =
\E\sup_{f_1\in\mathcal{F}_1}|\sum_{i=1}^n (f_1(X_i) - \E
f_1(X_i))|
\end{displaymath}
and
\begin{displaymath}
B = \E\sup_{f_2\in\mathcal{F}_2}|\sum_{i=1}^n f_2(X_i)|,
\end{displaymath}
we get by (\ref{firsttriangle}) and (\ref{etriangle}),
\begin{align}\label{aux1}
\p(&Z \ge (1 + \eta)\E Z +
t)\nonumber \\
 \le& \p(\sup_{f_1\in\mathcal{F}_1}|\sum_{i=1}^n (f_1(X_i) - \E f_1(X_i))| \ge (1
+ \eta)\E Z + (1-\varepsilon)t) \nonumber \\
&+ \p(\sup_{f_2\in\mathcal{F}_2}|\sum_{i=1}^n (f_2(X_i) - \E
f_2(X_i))| \ge \varepsilon
t)\nonumber\\
 \le& \p(\sup_{f_1\in\mathcal{F}_1}|\sum_{i=1}^n (f_1(X_i) - \E f_1(X_i))| \ge (1
+ \eta)A - 4 B + (1-\varepsilon)t) \nonumber \\
&+\p(\sup_{f_2\in\mathcal{F}_2}|\sum_{i=1}^n (f_2(X_i) - \E
f_2(X_i))| \ge \varepsilon t)
\end{align}
and similarly by (\ref{secondtriangle}) and (\ref{etriangle}),
\begin{align}\label{aux2}
\p(&Z \le (1-\eta)\E Z - t)  \nonumber\\
\le & \p(\sup_{f_1\in\mathcal{F}_1}|\sum_{i=1}^n (f_1(X_i) - \E
f_1(X_i))|  \le (1- \eta)\E Z -
(1-\varepsilon)t )\nonumber \\
&+ \p(\sup_{f_2\in\mathcal{F}_2}|\sum_{i=1}^n (f_2(X_i) - \E
f_2(X_i))| \ge \varepsilon t) \nonumber\\
\le& \p(\sup_{f_1\in\mathcal{F}_1}|\sum_{i=1}^n (f_1(X_i) - \E
f_1(X_i))|   \le (1- \eta)A -
(1-\varepsilon)t + 2 B)\nonumber \\
&+ \p(\sup_{f_2\in\mathcal{F}_2}|\sum_{i=1}^n (f_2(X_i) - \E
f_2(X_i))| \ge \varepsilon t).
\end{align}

We would like to choose a truncation level $\rho$ in a way, which
would allow to bound the first summands on the right-hand sides of
(\ref{aux1}) and (\ref{aux2}) with Lemma \ref{CLT_type} and the
other summands with Theorem \ref{Tal_psi}.

To this end let us set

\begin{align}\label{aux5}
\rho = 8 \E\max_{1\le i \le n}\sup_{f\in\mathcal{F}}| f(X_i)|\le
K_\alpha \Big\|\max_{1\le i\le n}\sup_{f\in \mathcal{F}}|
f(X_i)|\Big\|_{\psi_\alpha}.
\end{align}
Let us notice that by the Chebyshev inequality and the definition
of the class $\mathcal{F}_2$, we have

\begin{align*}
\p(\max_{k\le n} \sup_{f_2\in \mathcal{F}_2}|\sum_{i=1}^k
f_2(X_i)|
> 0) \le \p(\max_i\sup_f |f(X_i)| > \rho) \le 1/8
\end{align*}
and thus by the Hoffmann-J{\o}rgensen inequality (see e.g. \cite{LT},
Chapter 6, Proposition 6.8., inequality (6.8)), we obtain
\begin{align}\label{aux3}
B = \E\sup_{f_2 \in \mathcal{F}_2} |\sum_{i=1}^n f_2(X_i)| \le
8\E\max_{1\le i\le n}\sup_{f\in\mathcal{F}}|f(X_i)|.
\end{align}
In consequence
\begin{align*}
\E \sup_{f_2 \in\mathcal{F}_2} |\sum_{i=1}^n (f_2(X_i) - \E
f_2(X_i))|
&\le 16\E\max_{1\le i\le n}\sup_{f\in\mathcal{F}} |f(X_i)|\\
&\le K_\alpha\Big\|\max_{1\le i\le n}
\sup_{f\in\mathcal{F}}|f(X_i)| \Big\|_{\psi_\alpha}.
\end{align*}
We also have
\begin{align*}
\Big\|\max_{1\le i\le n} \sup_{f_2\in\mathcal{F}_2}|f_2(X_i) &- \E
f_2(X_i)| \Big\|_{\psi_\alpha} \\
&\le K_\alpha\Big\|\max_{1\le i\le n} \sup_{f_2\in\mathcal{F}_2}|f_2(X_i)|
\Big\|_{\psi_\alpha}  +
K_\alpha\Big\|\E \max_{1\le i\le n} \sup_{f_2\in\mathcal{F}_2}|f_2(X_i)|\Big\|_{\psi_\alpha} \\
&\le K_\alpha\Big\|\max_{1\le i\le n}
\sup_{f_2\in\mathcal{F}_2}|f_2(X_i)|\Big\|_{\psi_\alpha} \\
&\le K_\alpha\Big\|\max_{1\le i\le n}
\sup_{f\in\mathcal{F}}|f(X_i)|\Big\|_{\psi_\alpha}
\end{align*}
(recall that with our definitions, for $\alpha < 1$,
$\|\cdot\|_{\psi_\alpha}$ is a quasi-norm, which explains the
presence of the constant $K_\alpha$ in the first inequality).
Thus, by Theorem \ref{Tal_psi}, we obtain
\begin{displaymath}
\Big\| \sup_{f_2\in\mathcal{F}_2}|\sum_{i=1}^n (f_2(X_i) - \E
f_2(X_i))|\Big\|_{\psi_\alpha} \le K_\alpha \Big\|\max_{1\le i\le
n}\sup_{f\in\mathcal{F}} |f(X_i)| \Big\|_{\psi_\alpha},
\end{displaymath}
which implies
\begin{align}\label{aux4}
\p(\sup_{f_2\in\mathcal{F}_2}&|\sum_{i=1}^n f_2(X_i) - \E
f_2(X_i)|
\ge \varepsilon t) \nonumber\\
&\le 2\exp\Big(-\Big(\frac{\varepsilon t}{K_\alpha\|\max_{1\le
i\le
n}\sup_{f\in\mathcal{F}}|f(X_i)|\|_{\psi_\alpha}}\Big)^{\alpha}\Big).
\end{align}
Let us now choose $\varepsilon < 1/10$ and such that
\begin{align}\label{epsilon_choice}
(1-5\varepsilon)^{-2}(1+ \delta/2) \le (1 + \delta).
\end{align}

Since $\varepsilon$ is a function of $\delta$, in view of
(\ref{aux5}) and (\ref{aux3}), we can choose the constant
$K(\alpha,\eta,\delta)$ in (\ref{thewlogrestriction}) to be large
enough, to assure that
\begin{displaymath}
B \le 8\E\max_{1\le i\le n}\sup_{f\in\mathcal{F}}|f(X_i)| \le
\varepsilon t.
\end{displaymath}

Notice moreover, that for every $f\in \mathcal{F}$, we have $\E
(f_1(X_i) - \E f_1(X_i))^2 \le \E f_1(X_i)^2 \le \E f(X_i)^2$.

Thus, using inequalities (\ref{aux1}), (\ref{aux2}), (\ref{aux4})
and Lemma \ref{CLT_type} (applied for $\eta$ and $\delta/2$), we
obtain
\begin{align*}
\p(&Z \ge (1 + \eta)\E Z +
t), \quad \p(Z \le (1-\eta)\E Z - t) \nonumber\\
\le &
\exp\Big(-\frac{t^2(1-5\varepsilon)^2}{2(1+\delta/2)\sigma^2}\Big)
+
\exp\Big(-\frac{(1-5\varepsilon)t}{K(\eta,\delta) \rho}\Big) \\
&+ 2\exp\Big(-\Big(\frac{\varepsilon t}{K_\alpha\|\max_{1\le i\le
n}\sup_{f\in\mathcal{F}}|f(X_i)|\|_{\psi_\alpha}}\Big)^{\alpha}\Big).
\end{align*}

Since $\varepsilon < 1/10$, using (\ref{aux5}) one can see that
for $t$ satisfying (\ref{thewlogrestriction}) with
$K(\alpha,\eta,\delta)$ large enough,  we have
\begin{align*}
\exp&\Big(-\frac{(1-5\varepsilon)t}{K(\eta,\delta) \rho}\Big),
\exp\Big(-\Big(\frac{\varepsilon t}{K_\alpha\|\max_{1\le i\le
n}\sup_{f\in\mathcal{F}}|f(X_i)|\|_{\psi_\alpha}}\Big)^{\alpha}\Big)\\
&\le
\exp\Big(-\Big(\frac{t}{\tilde{C}(\alpha,\eta,\delta)\|\max_{1\le
i\le
n}\sup_{f\in\mathcal{F}}|f(X_i)|\|_{\psi_\alpha}}\Big)^{\alpha}\Big)
\end{align*}
(note that the above inequality holds for all $t$ if $\alpha =
1$).

Therefore, for such $t$,
\begin{align*}
\p(&Z \ge (1 + \eta)\E Z +
t), \quad \p(Z \le (1-\eta)\E Z - t) \nonumber\\
\le &
\exp\Big(-\frac{t^2(1-5\varepsilon)^2}{2(1+\delta/2)\sigma^2}\Big)\\
&+
3\exp\Big(-\Big(\frac{t}{\tilde{C}(\alpha,\eta,\delta)\|\max_{1\le
i\le n}\sup_{f\in
\mathcal{F}}|f(X_i)|\|_{\psi_\alpha}}\Big)^{\alpha}\Big).
\end{align*}
To finish the proof it is now enough to use
(\ref{epsilon_choice}).

\end{proof}

\paragraph{Remark} We would like to point out that the use of the Hoffman-J{\o}rgensen inequality in similar context is well known.
Such applications appeared in the proof of the aforementioned moment estimates for empirical processes by Gin\'{e}, Lata{\l}a, Zinn \cite{GLZ},
in the proof of Theorem \ref{Tal_psi} and recently in the proof of Fuk-Nagaev type inequalities for empirical processes used by Einmahl and Li to investigate generalized laws of the iterated logarithm for Banach space valued variables \cite{EL}.

As for using Theorem \ref{Tal_psi} to control the remainder after
truncating the original random variables, it was recently used in
a somewhat similar way by Mendelson and Tomczak-Jaegermann (see
\cite{MeT}).
\subsection{A counterexample \label{counter1}}

We will now present a simple example, showing that in Theorem
\ref{tail_estimate} one cannot replace $\|\sup_f \max_i
|f(X_i)|\|_{\psi_\alpha}$ with $\max_i\|\sup_f
|f(X_i)|\|_{\psi_\alpha}$. With such a modification, the
inequality fails to be true even in the real valued case, i.e.
when $\mathcal{F}$ is a singleton. For simplicity we will consider
only the case $\alpha = 1$.

Consider a sequence $Y_1,Y_2,\ldots,$ of i.i.d. real random
variables, such that $\p(Y_i = r) = e^{-r} = 1 - \p(Y_i = 0)$. Let
$\varepsilon_1,\varepsilon_2, \ldots,$ be a Rademacher sequence,
independent from $(Y_i)_i$. Define finally $X_i = \varepsilon_i
Y_i$. We have
\begin{displaymath}
\E e^{|X_i|} = e^{r}e^{-r} + (1 - e^{-r}) \le 2,
\end{displaymath}
so $\|X_i\|_{\psi_1} \le 1$. Moreover
\begin{displaymath}
\E |X_i|^2 = r^2e^{-r}.
\end{displaymath}

Assume now that we have for all $n,r \in \N$ and $t \ge 0$,
\begin{displaymath}
\p\Big(\Big|\sum_{i=1}^n X_i\Big| \ge K(\sqrt{nt}\|X_1\|_2 +
t\|X_1\|_{\psi_1})\Big) \le Ke^{-t},
\end{displaymath}
where $K$ is an absolute constant (which would hold if the
corresponding version of Theorem \ref{tail_estimate} was true).

For sufficiently large $r$, the above inequality applied with $n
\simeq e^r r^{-2}$ and $t \simeq r$, implies that
\begin{displaymath}
\p\Big(\Big|\sum_{i=1}^n X_i\Big| \ge r\Big) \le Ke^{-r/K}.
\end{displaymath}

On the other hand, by Levy's inequality, we have
\begin{displaymath}
2 \p\Big(\Big|\sum_{i=1}^n X_i\Big| \ge r\Big) \ge \p(\max_{i\le
n}|X_i| \ge r) \ge \frac{1}{2}\min(n\p(|X_1| \ge r),1) \ge
\frac{1}{2}r^{-2},
\end{displaymath}
which gives a contradiction for large $r$.

\paragraph{Remark} A small modification of the above argument shows that one cannot hope for an inequality
\begin{displaymath}
\p\Big(Z \ge K (\E Z + \sqrt{t}\sigma + t[\log^\beta
n]\max_i\|\sup_{f\in\mathcal{F}} |f(X_i)|\|_{\psi_1}\Big) \le K
e^{-t}
\end{displaymath}
with $\beta < 1$. For $\beta = 1$, this inequality follows from
Theorem \ref{tail_estimate} via Pisier's inequality \cite{P},
\begin{align}\label{Pisier}
\Big\|\max_{i\le n}|Y_i|\Big\|_{\psi_\alpha} \le K_\alpha \max_{i\le n} \|Y_i\|_{\psi_\alpha}\log^{1/\alpha} n
\end{align}
for independent real variables $Y_1$, \ldots, $Y_n$.

\section{Applications to Markov chains \label{Markov}}
We will now turn to the other class of inequalities we are
interested in. We are again concerned with random variables of
the form
\begin{displaymath}
Z = \sup_{f\in \mathcal{F}} |\sum_{i=1}^n f(X_i)|,
\end{displaymath}
but this time we assume that the class $\mathcal{F}$ is uniformly
bounded and we drop the assumption on the independence of the
underlying sequence $X_1,\ldots,X_n$. To be more precise, we will
assume that $X_1,\ldots,X_n$ form a Markov chain, satisfying some
additional conditions, which are rather classical in the Markov
chain or Markov Chain Monte Carlo literature.

The organization of this part is as follows. First, before stating
the main results, we will present all the structural assumptions
we will impose on the chain. At the same time we will introduce
some notation, which will be used in the sequel. Next, we present
our results (Theorems \ref{sums_Markov} and
\ref{empirical_Markov}) followed by the proof (which is quite
straightforward but technical) and a discussion of the optimality
(Section \ref{section_counterexample}). At the end, in Section
\ref{boundeddifference}, we will also present a bounded
differences type inequality for Markov chains.

\subsection{Assumptions on the Markov chain} Let $X_1,X_2,\ldots$
be a homogeneous Markov chain on $\mathcal{S}$, with transition
kernel $P = P(x,A)$, satisfying the so called \emph{minorization
condition}, stated below.

\paragraph{Minorization condition}
We assume that there exist positive $m \in \N$, $\delta > 0$, a
set $C \in \mathcal{B}$ (,,small set'') and a probability measure
$\nu$ on $\mathcal{S}$ for which
\begin{align}\label{minorization}
\forall_{x \in C} \; \forall_{A \in\mathcal{B}} \; P^m(x,A) \ge
\delta\nu(A)
\end{align}
and
\begin{align}\label{access}
\forall_{x \in \mathcal{S}} \exists_n \; P^{nm}(x,C) > 0,
\end{align}
where $P^i(\cdot,\cdot)$ is the transition kernel for the chain after $i$ steps.

One can show that in such a situation if the chain admits an
invariant measure $\pi$, then this measure is unique and satisfies
 $\pi(C)
> 0$ (see \cite{MT}). Moreover, under some conditions on the initial distribution
$\xi$, it can be extended to a new (so called \emph{split}) chain
$(\tilde{X}_n,R_n) \in \mathcal{S}\times\{0,1\}$, satisfying the
following properties.

\paragraph{Properties of the split chain}
\begin{itemize}
\item[(P1)] $(\tilde{X}_n)_n$ is again a Markov chain with
transition kernel $P$ and initial distribution $\xi$ (hence for
our purposes of estimating the tail probabilities we may and will
identify $X_n$ and $\tilde{X}_n$ ),

\item[(P2)] if we define $T_1 = \inf\{n>0\colon R_{nm} = 1\}$,
\begin{displaymath}
T_{i+1} = \inf\{n>0\colon R_{(T_1+\ldots + T_i+n)m} = 1\},
\end{displaymath}
then $T_1,T_2,\ldots,$ are well defined, independent, moreover
$T_2,T_3,\ldots$ are i.i.d.,

\item[(P3)] if we define $S_i = T_1+\ldots+T_i$, then the
,,blocks''
\begin{align*}
Y_0 &= (X_1,\ldots,X_{mT_1 + m -1}),\\
Y_i & = (X_{m(S_i + 1)},\ldots,X_{mS_{i+1}+m-1}),\quad i > 0,
\end{align*}
form a one-dependent sequence (i.e. for all $i$, $\sigma((Y_j)_{j<
i})$ and $\sigma((Y_j)_{j> i})$ are independent). Moreover, the
sequence $Y_1,Y_2,\ldots$ is stationary. If $m=1$, then the
variables $Y_0,Y_1,\ldots$ are independent.

 In consequence, for $f\colon \mathcal{S} \to \R$, the variables
\begin{displaymath}
Z_i = Z_i(f) = \sum_{i=m(S_i +1)}^{mS_{i+1} + m -1} f(X_i), \; i
\ge 1,
\end{displaymath}
constitute a one-dependent stationary sequence (an i.i.d. sequence
if $m=1$). Additionally, if $f$ is $\pi$-integrable (recall that $\pi$ is the unique stationary measure for the chain), then
\begin{align}\label{integration_formula}
\E Z_i = \delta^{-1}\pi(C)^{-1}m\int f d\pi.
\end{align}

\item[(P4)] the distribution of $T_1$ depends only on
$\xi,P,C,\delta,\nu$, whereas the law of $T_2$ only on
$P,C,\delta$ and $\nu$.
\end{itemize}

We refrain from specifying the construction of this new chain in
full generality as well as conditions under which
(\ref{minorization}) and (\ref{access}) hold and refer the reader
to the classical monograph \cite{MT} or a survey article \cite{RR}
for a complete exposition. Here, we will only sketch the
construction for $m=1$, to give its ,,flavour''. Informally
speaking, at each step $i$, if we have $X_i = x$ and $x \notin C$,
we generate the next value of the chain, according to the measure
$P(x,\cdot)$. If $x \in C$, then we toss a coin with probability
of success equal to $\delta$. In the case of success ($R_i = 1$),
we draw the next sample according to the measure $\nu$, otherwise
($R_i=0)$, according to
\begin{displaymath}
\frac{P(x,\cdot) - \delta\nu(\cdot)}{1 - \delta}.
\end{displaymath}

When $R_i = 1$, one usually says that the chain
\emph{regenerates}, as the distribution in the next step (for
$m=1$, after $m$ steps in general) is again $\nu$.

Let us remark that for a recurrent chain on a countable state
space, admitting a stationary distribution, the Minorization
condition is always satisfied with $m = 1$ and $\delta = 1$ (for
$C$ we can take $\{x\}$, where $x$ is an arbitrary element of the
state space). Also the construction of the split chain becomes
trivial.

\paragraph{}
Before we proceed, let us present a general idea, our approach is based on. To derive our estimates we will need two types of assumptions.
\paragraph{Regeneration assumption.}
We will work under the assumption that the chain admits a
representation as above (Properties (P1) to (P4)). We will not
however take advantage of the explicit construction. Instead we
will use the properties stated in points above. A similar approach
is quite common in the literature.

\paragraph{Assumption of the exponential integrability of the regeneration time.}
To derive concentration of measure inequalities, we will also
assume that $\|T_1\|_{\psi_1} < \infty$ and $\|T_2\|_{\psi_1}<
\infty$. At the end of the article we will present examples for
which this assumption is satisfied and relate obtained
inequalities to known results.

\paragraph{}
The regenerative properties of the chain allow us to decompose the
chain into one-dependent (independent if $m=1$) blocks of random
length, making it possible to reduce the analysis of the chain to
sums of independent random variables (this approach is by now
classical, it has been successfully used in the analysis of limit
theorems for Markov chains, see \cite{MT}). Since we are
interested in non-asymptotic estimates of exponential type, we
have to impose some additional conditions on the regeneration
time, which would give us control over the random length of
one-dependent blocks. This is the reason for introducing the
assumption of the exponential integrability which (after some
technical steps) allows us to apply the inequalities for unbounded
empirical processes, presented in Section \ref{independent}.

\subsection{Main results concerning Markov chains\label{subsectionMarkov}}
\paragraph{}
Having established all the notation, we are ready to state our
main results on Markov chains.

As announced in the introduction, our results depend on the
parameter $m$ in the Minorization condition. If $m=1$ we are able
to obtain tail inequalities for empirical processes (Theorem
\ref{empirical_Markov}), whereas for $m > 1$ we have to restrict
to linear statistics of Markov chains (Theorem \ref{sums_Markov}),
which formally corresponds to empirical processes indexed by a
singleton. The variables $T_i$ and $Z_1$ appearing in the theorems were
defined in the previous section (see the properties (P2) and (P3) of the
split chain).

\begin{theorem}\label{sums_Markov}
 Let $X_1,X_2,\ldots$ be a Markov chain with values in $\mathcal{S}$, satisfying the \textbf{Minorization condition} and admitting a (unique)
stationary distribution $\pi$. Assume also that $\|T_1\|_{\psi_1},\|T_2\|_{\psi_1} \le \tau$. Consider
 a function
$f\colon \mathcal{S} \to \R$, such that $\|f\|_\infty \le a$ and
$\E_\pi f = 0$. Define also the random variable
\begin{align*}
Z = \sum_{i=1}^n f(X_i).
\end{align*}
Then for all $t > 0$,
\begin{align}\label{sums_Markov_eq}
\p\Big(|Z| > t \Big) \le&
K\exp\Big(-\frac{1}{K}\min\Big(\frac{t^2}{n (m\E
T_2)^{-1}\Var{Z_1}},\frac{t}{\tau^2am\log n}\Big)\Big).
\end{align}
\end{theorem}

\begin{theorem} \label{empirical_Markov} Let $X_1,X_2,\ldots$ be a Markov chain with values in $\mathcal{S}$, satisfying
the \textbf{Minorization condition} with $m = 1$ and admitting a (unique) stationary distribution $\pi$.
Assume also that $\|T_1\|_{\psi_1},\|T_2\|_{\psi_1} \le \tau$.
 Consider moreover a countable class $\mathcal{F}$ of
measurable functions $f \colon \mathcal{S} \to \R$, such that
$\|f\|_\infty \le a$ and $\E_\pi f = 0$. Define the random
variable
\begin{displaymath}
Z = \sup_{f\in\mathcal{F}}\Big|\sum_{i=1}^n f(X_i) \Big|
\end{displaymath}
and the ''asymptotic weak variance''
\begin{displaymath}
\sigma^2 = \sup_{f\in \mathcal{F}} \Var{Z_1(f)}/\E T_2.
\end{displaymath}
Then, for all $t \ge  1$,
\begin{displaymath}
\p\Big(Z \ge K\E Z + t\Big) \le
K\exp\Big(-\frac{1}{K}\min\Big(\frac{t^2}{n\sigma^2},\frac{t}{\tau^3(\E
T_2)^{-1} a \log n}\Big)\Big).
\end{displaymath}
\end{theorem}

\paragraph{Remarks}
\begin{enumerate}
\item As it was mentioned in the previous section, chains satisfying the Minorization condition
admit at most one stationary measure.

\item In Theorem \ref{empirical_Markov}, the dependence on the
chain is worse that in Theorem \ref{sums_Markov}, i.e. we have
$\tau^3(\E T_2)^{-1}$ instead of $\tau^2$ in the denominator. It
is a result of just one step in the argument we present below,
however at the moment we do not know how to improve this
dependence (or extend the result to $m>1$).

\item Another remark we would like to make is related to the limit behaviour of the Markov chain.
Let us notice that the asymptotic variance (the variance in the CLT) for $n^{-1/2}(f(X_1)+\ldots+f(X_n))$ equals
\begin{displaymath}
m^{-1}(\E T_2)^{-1} (\Var Z_1 + \E Z_1Z_2),
\end{displaymath}
which for $m = 1$ reduces to
\begin{displaymath}
(\E T_2)^{-1} \Var Z_1
\end{displaymath}
(we again refer the reader to \cite{MT}, Chapter 17 for details).
Thus, for $m=1$ our estimates reflect the asymptotic behaviour of the variable $Z$.
\end{enumerate}

\paragraph{}
Let us now pass to the proofs of the above theorems. For a function $f \colon \mathcal{S} \to \R$ let us define
\begin{displaymath}
Z_0 = \sum_{i=1}^{(mT_1+m-1)\wedge n} f(X_i)
\end{displaymath}
and recall the variables $S_i$ and
\begin{displaymath}
Z_i = Z_i(f) = \sum_{i=m(S_i +1)}^{mS_{i+1} + m -1} f(X_i), \; i
\ge 1,
\end{displaymath}
defined in the previous section (see property (P3) of the split
chain). Recall also that $Z_i$'s form a one-dependent stationary
sequence for $m > 1$ and an i.i.d. sequence for $m=1$.

Using this notation, we have
\begin{align}\label{block_decomposition}
 f(X_1)+\ldots+f(X_n) = Z_0 + \ldots + Z_{N} + \sum_{i=
(S_{N+1} + 1)m}^n f(X_i),
\end{align}
with
\begin{align}\label{Definition_of_n}
N = \sup\{i\in \N\colon m S_{i+1} + m - 1\le n\},
\end{align}
 where $\sup\emptyset = 0$ (note that $N$ is a random variable). Thus
$Z_0$ represents the sum up to the first regeneration time, then
$Z_1,\ldots,Z_{N}$ are identically distributed blocks between
consecutive regeneration times, included in the interval $[1,n]$,
finally the last term corresponds to the initial segment of the
last block. The sum $Z_1+\ldots+Z_{N}$ is empty if up to time $n$,
there has not been any regeneration (i.e. $mT_1 + m -1 > n$) or
there has been only one regeneration ($mT_1 + m -1 \le n$ and
$m(T_1+T_2) + m - 1 > n$). The last sum on the right hand side is
empty if there has been no regeneration or the last 'full' block
ends with $n$.

We will first bound the initial and the last summand in the
decomposition (\ref{block_decomposition}). To achieve this we will
not need the assumptions that $f$ is centered with respect to the
stationary distribution $\pi$. In consequence the same bound may
be applied to proofs of both Theorem \ref{sums_Markov} and Theorem
\ref{empirical_Markov}.

\begin{lemma}\label{initial_block}
If $\|T_1\|_{\psi_1}\le \tau$ and $\|f\|_\infty \le a$, then for
all $t \ge 0$,
\begin{displaymath}
\p(|Z_0| \ge t) \le 2\exp\Big(\frac{-t}{2am\tau}\Big).
\end{displaymath}
\end{lemma}

\begin{proof}
We have $|Z_0|  \le 2aT_1m$, so by the remark after Definition
\ref{Orlicz_norm},
\begin{displaymath}
\p(|Z_0| \ge t) \le \p(T_1 \ge t/2am) \le 2\exp\Big(\frac{-t}{2am\tau}\Big).
\end{displaymath}
\end{proof}

The next lemma provides a similar bound for the last summand on
the right hand side of (\ref{block_decomposition}). It is a little
bit more complicated, since it involves additional dependence on
the random variable $N$.

\begin{lemma}\label{thelastblock}
If $\|T_1\|_{\psi_1}, \|T_2\|_{\psi_1} \le \tau$, then for all $t \ge 0$,
\begin{displaymath}
\p(n - m(S_{N+1}+1) + 1 > t) \le K\exp\Big(\frac{-t}{Km \tau\log\tau}\Big).
\end{displaymath}
In consequence, if $\|f\|_\infty
\le a$, then
\begin{displaymath}
\p\Big(\Big|\sum_{i= (S_{N+1} + 1)m}^n f(X_i)\Big| > t\Big) \le
K\exp\Big(\frac{-t}{Kam\tau\log\tau}\Big).
\end{displaymath}
\end{lemma}
\begin{proof}
Let us consider the variable $M_n = n - m(S_{N+1}+1) + 1$. If $M_n
> t$ then
\begin{displaymath}
S_{N+1} < \frac{n-t + 1}{m} - 1 < \Big\lfloor\frac{n-t + 1}{m}\Big\rfloor.
\end{displaymath}

Therefore
\begin{align*}
\p(M_n > t) &\le \sum_{k< \frac{n-t + 1}{m} - 1}\p(S_{N+1} = k)\\
&= \sum_{k<\frac{n-t + 1}{m} - 1}\sum_{l=1}^k \p(S_l = k\; \& {N+1} = l) \\
&= \sum_{k<\frac{n-t + 1}{m} - 1}\sum_{l=1}^k \p(S_l = k\; \& \;m(k + T_{l+1}) + m - 1 > n)\\
&=\sum_{k<\frac{n-t + 1}{m} - 1}\sum_{l=1}^k\p(S_l = k) \p(T_2 > \frac{n+1}{m}-1-k) \\
&\le\sum_{k=1}^{\lfloor\frac{n-t + 1}{m}\rfloor-1}\p(T_2 > \frac{n+1}{m}-1-k)\\
&\le \sum_{k=1}^{\lfloor\frac{n-t + 1}{m}\rfloor-1} 2\exp\Big(\frac{1}{\tau}\Big(k+1-\frac{n+1}{m}\Big)\Big) \\
&\le 2\exp\Big(\tau^{-1}\Big(1 - \frac{n+1}{m}\Big)\Big)\exp(1/\tau)\frac{\exp\Big(\tau^{-1}\Big(\lfloor\frac{n-t + 1}{m}\rfloor - 1\Big)\Big)}{\exp(1/\tau) - 1}\\
&\le K\tau \exp\Big(\frac{-t}{m\tau}\Big),
\end{align*}
where the first equality follows from the fact that $S_{N+1} \ge
N+1$, the second from the definition of $N$, the third from the
fact that $T_1,T_2,\ldots$ are independent and $T_2,T_3,\ldots$
are i.i.d., finally the second inequality from the fact that
$S_i\neq S_j$ for $i\neq j$ (see the properties (P2) and (P3) of the split chain).

Let us notice that if $t > 2m\tau\log \tau$, then
\begin{displaymath}
\tau\exp\Big(\frac{-t}{m\tau}\Big) \le \exp\Big(\frac{-t}{2m\tau}\Big) \le \exp\Big(\frac{-t}{Km\tau\log \tau}\Big),
\end{displaymath}
where in the last inequality we have used the fact that $\tau > c$ for some universal constant $c > 1$.

On the other hand, if $t < 2m\tau\log\tau$, then
\begin{displaymath}
1 \le e\cdot\exp\Big(\frac{-t}{2m \tau\log\tau}\Big).
\end{displaymath}
Therefore we obtain for $t \ge 0$,
\begin{align*}
\p(M_n > t) \le K\exp\Big(\frac{-t}{Km \tau\log\tau}\Big),
\end{align*}
which proves the first part of the Lemma.
Now,
\begin{displaymath}
\p\Big(\Big|\sum_{i= (S_{N+1} + 1)m}^n f(X_i)\Big| > t\Big) \le \p(M_n > t/a)\le K\exp\Big(\frac{-t}{Kam\tau\log\tau}\Big)
\end{displaymath}
\end{proof}

Before we proceed with the proof of Theorem \ref{sums_Markov} and
Theorem \ref{empirical_Markov}, we would like to make some
additional comments regarding our approach. As already mentioned, thanks to the property
(P3) of the split chain, we may apply to $Z_1+\ldots+Z_N$ the
inequalities for sums of independent random variables obtained in
Section \ref{independent} (since for $m > 1$ we have only
one-dependence, we will split the sum, treating even and odd
indices separately). The number of summands is random, but clearly
not larger than $n$. Since the variables $Z_i$ are
equidistributed, we can reduce this random sum to a deterministic
one by applying the following maximal inequality by
Montgomery-Smith \cite{MS}.

\begin{lemma}\label{Montgomery_Smith}
Let $Y_1,\ldots,Y_n$ be i.i.d. Banach space valued random
variables. Then for some universal constant $K$ and every $t > 0$,
\begin{displaymath}
\p\Big(\max_{k\le n} \Big\|\sum_{i=1}^k Y_i \Big\| > t \Big) \le K
\p\Big(\Big\|\sum_{i=1}^n Y_i \Big\| > t/K \Big).
\end{displaymath}
\end{lemma}

\paragraph{Remark} The use of regeneration methods makes our proof similar to
the proof of the CLT for Markov chains. In this context, the above
lemma can be viewed as a counterpart of the Anscombe theorem (they
are quite different statements but both are used to handle the
random number of summands).

One could now apply Lemma \ref{Montgomery_Smith} directly, using the fact that $N \le n$. Then however one would
not get the asymptotic variance in the exponent (see remark after Theorem \ref{empirical_Markov}). The form of this variance
is a consequence of the aforementioned Anscombe theorem and the fact that by the LLN we have (denoting $N = N_n$ to stress the
dependence on $n$)
\begin{align}\label{LLN}
\lim_{n\to \infty} \frac{N_n}{n} = \frac{1}{m\E T_2}\quad {\rm
a.s.}
\end{align}

Therefore to obtain an inequality which at least up to universal constants (and for $m=1$) reflects the limiting behaviour of the variable $Z$, we
will need a quantitative version of (\ref{LLN}) given in the following lemma.

\begin{lemma}\label{estimates_on_N}
If $\|T_1\|_{\psi_1}, \|T_2\|_{\psi_1} \le \tau$, then
\begin{displaymath}
\p( N >  \lfloor 3n/(m\E T_2)\rfloor ) \le
K\exp\Big(-\frac{1}{K}\frac{n\E T_2}{m\tau^2}\Big).
\end{displaymath}
\end{lemma}

To prove the above estimate, we will use the classical
Bernstein's inequality (actually its version for $\psi_1$
variables).
\begin{lemma}[Bernstein's $\psi_1$ inequality, see \cite{VW}, Lemma 2.2.11 and the subsequent remark]
 \label{Bernstei_pis_1}
If $Y_1,\ldots,Y_n$ are independent random variables such that $\E Y_i = 0$ and $\|Y\|_{\psi_1} \le \tau$, then for every $t > 0$,
\begin{align*}
\p\Big(\Big| \sum_{i=1}^n Y_i\Big| > t\Big) \le
2\exp\Big(-\frac{1}{K}\min\Big(\frac{t^2}{n\tau^2},
\frac{t}{\tau}\Big)\Big).
\end{align*}
\end{lemma}

\begin{proof}[Proof of Lemma \ref{estimates_on_N}]

Assume now that $n/(m\E T_2) \ge 1$. We have
\begin{align*}
\p( N > \lfloor 3&n/(m\E T_2)\rfloor) \le \p\Big( m (T_2+\ldots+ T_{\lfloor 3 n/(m\E T_2)\rfloor + 1})\le n\Big)\\
&\le \p\Big(\sum_{i=2}^{\lfloor 3n/(m\E T_2)\rfloor+1} (T_i - \E T_2) \le n/m - \lfloor 3n/(m\E T_2)\rfloor \E T_2\Big)\\
&\le \p\Big(\sum_{i=2}^{\lfloor 3n/(m\E T_2)\rfloor+1} (T_i - \E T_2) \le n/m - 3n/(2 m)\Big)\\
& = \p\Big(\sum_{i=2}^{\lfloor 3n/(m\E T_2)\rfloor+1} (T_i - \E T_2) \le -n/(2m)\Big).
\end{align*}
We have $\|T_2 - \E T_2\|_{\psi_1} \le 2 \|T_2\|_{\psi_1} \le 2\tau$, therefore Bernstein's inequality (Lemma \ref{Bernstei_pis_1}), gives
\begin{align*}
\p( N >  \lfloor 3&n/(m\E T_2)\rfloor ) \le 2\exp\Big(-\frac{1}{K}\min\Big(\frac{(n/2m)^2}{(3 n/(m\E T_2))\tau^2},\frac{n}{m\tau}\Big)\Big)\\
&\le2\exp\Big(-\frac{1}{K}\min\Big(\frac{n\E T_2}{m\tau^2},\frac{n}{m\tau}\Big)\Big) \nonumber\\
&= 2\exp\Big(-\frac{1}{K}\frac{n\E T_2}{m\tau^2}\Big),\nonumber
\end{align*}
where the equality follows from the fact that $\E T_2 \le \tau$.
If $n/(m\E T_2) < 1$, then also $n\E T_2/(m\tau^2) < 1$, thus
finally we have
\begin{displaymath}
\p( N >  \lfloor 3n/(m\E T_2)\rfloor ) \le
K\exp\Big(-\frac{1}{K}\frac{n\E T_2}{m\tau^2}\Big),
\end{displaymath}
which proves the lemma.
\end{proof}

We are now in position to prove Theorem \ref{sums_Markov}

\begin{proof}[Proof of Theorem \ref{sums_Markov}]
Let us notice that $|Z_i| \le a mT_{i+1}$, so for $i\ge 1$, $\|Z_i\|_{\psi_1}
\le am\|T_2\|_{\psi_1} \le am\tau$. Additionally, by (\ref{integration_formula}), for $i \ge 1$, $\E Z_i = 0$. Denote now $R = \lfloor
3n/(m\E T_2)\rfloor$. Lemma \ref{estimates_on_N}, Lemma
\ref{Montgomery_Smith} and Theorem \ref{tail_estimate} (with
$\alpha = 1$, combined with Pisier's estimate (\ref{Pisier})) give
\begin{align}\label{middle_segment}
\p\Big( |Z_1 &+ \ldots + Z_N| > 2t\Big) \\
\le& \p\Big( |Z_1 + \ldots + Z_N| > 2t \; \& N \;\le R \Big)
+ K\exp\Big(-\frac{1}{K}\frac{n\E T_2}{m\tau^2}\Big)\nonumber\\
\le&\p\Big( |Z_1 + Z_3 + \ldots + Z_{2\lfloor (N -1)/2\rfloor+1}| > t \; \& N \;\le R \Big)\nonumber\\
&\p\Big( |Z_2 + Z_4 + \ldots + Z_{2\lfloor N /2\rfloor}| > t \; \& N \;\le R \Big)+ K\exp\Big(-\frac{1}{K}\frac{n\E T_2}{m\tau^2}\Big)\nonumber\\
\le& \p\Big( \max_{k\le\lfloor (R-1)/2\rfloor} |Z_1+Z_3+\ldots+Z_{2k+1}| > t \Big) \nonumber\\
&+ \p\Big( \max_{k\le\lfloor R/2\rfloor} |Z_2+\ldots+Z_{2k}| > t \Big)+ K\exp\Big(-\frac{1}{K}\frac{n\E T_2}{m\tau^2}\Big)\nonumber\\
\le& K\p\Big( |Z_1+Z_3+ \ldots+Z_{2\lfloor(R-1)/2\rfloor+1}| > t/K \Big) \nonumber\\
&+ K\p\Big( |Z_2+Z_4+ \ldots+Z_{2\lfloor R/2\rfloor}| > t/K \Big)+ K\exp\Big(-\frac{1}{K}\frac{n\E T_2}{m\tau^2}\Big)\nonumber\\
\le& K\exp\Big(-\frac{1}{K}\min\Big(\frac{t^2}{n (m\E T_2)^{-1}\Var{Z_1}},\frac{t}{\log(3n/(m \E T_2))am\tau}\Big)\Big) \nonumber \\
&+ K\exp\Big(-\frac{1}{K}\frac{n\E T_2}{m\tau^2}\Big)\nonumber.
\end{align}
Combining the above estimate with (\ref{block_decomposition}),
Lemma \ref{initial_block} and Lemma \ref{thelastblock}, we obtain
\begin{align*}
\p\Big(|S_n| &> 4t \Big) \\
\le& K\exp\Big(-\frac{1}{K}\min\Big(\frac{t^2}{n (m\E T_2)^{-1}\Var{Z_1}},\frac{t}{\log(3n/(m \E T_2))am\tau}\Big)\Big)\\
&+ K\exp\Big(-\frac{1}{K}\frac{n\E T_2}{m\tau^2}\Big) +
2\exp\Big(\frac{-t}{2am\tau}\Big) +
K\exp\Big(\frac{-t}{Kam\tau\log\tau}\Big).
\end{align*}
For $t > na/4$, the left hand side of the above inequality is equal to $0$, therefore, using the fact that $\E T_2\ge 1,\tau > 1$, we obtain
(\ref{sums_Markov_eq}).
\end{proof}

The proof of Theorem \ref{empirical_Markov} is quite similar,
however it involves some additional technicalities related to the
presence of $\E Z$ in our estimates.

\begin{proof}[Proof of Theorem \ref{empirical_Markov}]

Let us first notice that similarly as in the real valued case,
we have
\begin{align*}
\p(\sup_f |Z_0(f)| \ge t) \le \p(T_1 \ge t/a) \le
2\exp\Big(\frac{-t}{a\tau}\Big),
\end{align*}
moreover Lemma \ref{thelastblock} (applied to
the function $x\mapsto \sup_{f\in \mathcal{F}} |f(x)|$) gives
\begin{align*}
\p\Big(\sup_f \Big|\sum_{i= S_{N+1} + 1}^n f(X_i)\Big| > t\Big)
\le K\exp\Big(\frac{-t}{Ka\tau\log\tau}\Big).
\end{align*}

One can also see that since we assume that $m = 1$, the splitting of
$Z_1+\ldots+Z_N$ into sums over even and odd indices is not
necessary (by the property (P3) of the split chain the summands
are independent). Using the fact that Lemma \ref{Montgomery_Smith}
is valid for Banach space valued variables, we can repeat the
argument from the proof of Theorem \ref{sums_Markov} and obtain
for $R= \lfloor 3n/\E T_2\rfloor$,
\begin{align*}
\p\Big(Z \ge K\E \sup_f \Big|\sum_{i=1}^R Z_i(f)\Big| + &t\Big)
\le
K\exp\Big(-\frac{1}{K}\min\Big(\frac{t^2}{n\sigma^2},\frac{t}{\tau^2
a \log n}\Big)\Big)\\
&\le
K\exp\Big(-\frac{1}{K}\min\Big(\frac{t^2}{n\sigma^2},\frac{t}{\tau^3(\E
T_2)^{-1} a \log n}\Big)\Big).
\end{align*}
Thus, Theorem \ref{empirical_Markov} will follow if we prove that
\begin{align}\label{toprove}
\E \sup_f \Big|\sum_{i=1}^R Z_i(f)\Big| \le
K\E\sup_f\Big|\sum_{i=1}^n f(X_i)\Big| + K\tau^3 a/\E T_2
\end{align}
(recall that $K$ may change from line to line).

From the triangle inequality, the fact that $Y_i=
(X_{S_i+1},\ldots,X_{S_{i+1}})$, $i\ge1$, are i.i.d. and Jensen's
inequality it follows that
\begin{align}\label{triange}
\E\sup_f\Big|\sum_{i=1}^R Z_i(f) \Big| &\le
12\E\sup_f\Big|\sum_{i=1}^{\lceil n/(4\E T_2)\rceil}
Z_i(f)\Big|\nonumber \\
&\le 12\E\sup_f\Big|\sum_{i=1}^{\lfloor n/(4\E T_2)\rfloor}
Z_i(f)\Big| + 12a\tau,
\end{align}
where in the last inequality we used the fact that
$\E\sup_f|Z_i(f)|\le \E a T_{i+1} \le a\tau$.

 We will split the integral on the right hand side
into two parts, depending on the size of the variable $N$. Let us
first consider the quantity
\begin{align*}
\E\sup_f\Big|\sum_{i=1}^{\lfloor n/(4\E T_2)\rfloor}
Z_i(f)\Big|\ind{N < \lfloor n/(4\E T_2)\rfloor}
\end{align*}

Assume that $n/(4\E T_2) \ge 1$. Then, using Bernstein's
inequality, we obtain
\begin{align*}
\p\Big(N &< \lfloor n/(4\E T_2)\rfloor\Big) =
\p\Big(\sum_{i=1}^{\lfloor n/(4\E T_2)\rfloor+1} T_i > n \Big)
\\
&\le \p(T_1 >n/2) + \p\Big(\sum_{i=2}^{\lfloor n/(4\E
T_2)\rfloor+1} (T_i -
\E T_2) > n/2 - \lfloor n/(4\E T_2)\rfloor\E T_2\Big)\\
&\le 2e^{-n/2\tau} + \p\Big(\sum_{i=2}^{\lfloor n/(4\E
T_2)\rfloor+1} (T_i - \E
T_2) > n/4\Big)\\
&\le 2e^{-n/2\tau} + 2\exp\Big(-\frac{1}{K}\min\Big(\frac{n^2\E
T_2}{n\tau^2},\frac{n}{\tau}\Big)\Big) \le Ke^{-n\E T_2/K\tau^2}.
\end{align*}
If $(n/4\E T_2) < 1$, the above estimate holds trivially.
Therefore
\begin{align}\label{first_part}
\E\sup_f\Big|&\sum_{i=1}^{\lfloor n/(4\E T_2)\rfloor}
Z_i(f)\Big|\ind{N < \lfloor n/(4\E T_2)\rfloor} \le
a\sum_{i=1}^{\lfloor n/(4\E T_2)\rfloor} \E T_{i+1}\ind{N <
\lfloor
n/(4\E T_2)\rfloor}\nonumber\\
&\le an\|T_2\|_2\sqrt{\p(N < \lfloor n/(4\E T_2)\rfloor)}\nonumber\\
& \le Ka \tau ne^{-n\E T_2/K\tau^2} \le K a\tau^3/\E T_2.
\end{align}

Now we will bound the remaining part i.e.

\begin{align*}
\E\sup_f\Big|\sum_{i=1}^{\lfloor n/(4\E T_2)\rfloor}
Z_i(f)\Big|\ind{N \ge \lfloor n/(4\E T_2)\rfloor}.
\end{align*}

Recall that $Y_0 = (X_1,\ldots,X_{T_1})$, $Y_i= (X_{S_i+1},\ldots,
X_{S_{i+1}})$ for $i \ge 1$ and consider a filtration
$(\mathcal{F}_i)_{i\ge 0}$ defined as
\begin{displaymath}
\mathcal{F}_i = \sigma(Y_0,\ldots,Y_i),
\end{displaymath}
where we regard the blocks $Y_i$ as random variables with values
in the disjoint union $\bigcup_{i=1}^\infty \mathcal{S}^i$, with the
natural $\sigma$-field, i.e. the $\sigma$-field generated by
$\bigcup_{i=1}^\infty \mathcal{B}^{\otimes i}$ (recall that
$\mathcal{B}$ denotes our $\sigma$-field of reference in
$\mathcal{S}$).

Let us further notice that $T_i$ is measurable with respect to
$\sigma(Y_{i-1})$ for $i \ge 1$. We have for $i \ge 1$,
\begin{displaymath}
\{N + 1 \le i\} = \{T_1+\ldots+T_{i+1} > n\} \in \mathcal{F}_i
\end{displaymath}
and $\{N +1 \le 0\} = \emptyset$, so $N+1$ is a stopping time with
respect to the filtration $\mathcal{F}_i$. Thus we have
\begin{align*}
\E\sup_f\Big|&\sum_{i=1}^{\lfloor n/(4\E T_2)\rfloor}
Z_i(f)\Big|\ind{N \ge \lfloor n/(4\E T_2)\rfloor} \\
&=  \E\sup_f\Big|\sum_{i=1}^{\lfloor n/(4\E T_2)\rfloor\wedge
(N+1)} Z_i(f)\Big|\ind{N +1> \lfloor n/(4\E
T_2)\rfloor}\\
&\le \E \Big(\E\Big[\sup_f\Big|\sum_{i=1}^{N+1}
Z_i(f)\Big||\mathcal{F}_{\lfloor n/(4\E T_2)\rfloor\wedge
(N+1)}\Big]\ind{N
+1> \lfloor n/(4\E T_2)\rfloor}\Big)\\
&= \E \sup_f\Big|\sum_{i=1}^{N+1} Z_i(f)\Big|\ind{N +1>
\lfloor n/(4\E T_2)\rfloor}\\
&\le a\tau + \E \sup_f\Big|\sum_{i=0}^{N+1} Z_i(f)\Big|\ind{N +1>
\lfloor n/(4\E T_2)\rfloor}\\
&\le a\tau + \E \sup_f\Big|\sum_{i=1}^n f(X_i)\Big| +
\E\sup_f\Big|\sum_{i=n+1}^{S_{N+2}} f(X_i)\Big|,
\end{align*}
where in the first inequality we used Doob's optional sampling
theorem together with the fact that $\sup_f|\sum_{i=1}^n Z_i(f)|$
is a submartingale with respect to $(\mathcal{F}_i)$ (notice that
$Z_i(f)$ is measurable with respect to $\sigma(Y_i)$ for $i \in
\N$ and $f\in\mathcal{F}$). The second equality follows from the
fact that $\{N+1 > \lfloor n/(4\E T_2)\rfloor\} \in
\mathcal{F}_{\lfloor n/(4\E T_2)\rfloor \wedge (N+1)}$. Indeed for
$i \ge \lfloor n/(4\E T_2)\rfloor$, we have
\begin{align*}
&\{N+1 > \lfloor n/(4\E T_2)\rfloor \;\&\; \lfloor n/(4\E
T_2)\rfloor\wedge (N+1) \le i\} \\
&= \{N+1> \lfloor n/(4\E T_2)\rfloor\}  \in \mathcal{F}_{\lfloor
n/(4\E T_2)\rfloor} \subseteq \mathcal{F}_i,
\end{align*}
whereas for $i < \lfloor n/(4\E T_2)\rfloor$, this set is empty.

\paragraph{}
 Now, combining the above estimate with (\ref{triange}) and
(\ref{first_part}) and taking into account the inequality $\tau
\ge \E T_2 \ge  1$, it is easy to see that to finish the proof of
(\ref{toprove}) it is enough to show that
\begin{align}\label{tofinish}
\E\sup_f\Big|\sum_{i=n+1}^{S_{N+2}} f(X_i)\Big| \le Ka\tau^3/\E
T_2
\end{align}

This in turn will follow if we prove that $\E(S_{N+2} - n) \le
K\tau^3/\E T_2$.

Recall now the first part of Lemma \ref{thelastblock}, stating under our assumptions ($m=1$) that
$\p(n - S_{N+1} > t) \le K\exp(-t/K\tau\log \tau)$ for $t\ge 0$.
We have
\begin{align*}
\p(S_{N+2}& - n > t)  \\
\le&\p(n - S_{N+1} > t) + \p(n - S_{N+1} \le
t \;\&\; S_{N+2}-n>t \; \&N >0)\\
&+\p(S_{N+2} - n > t\;\&\; N=0)\\
\le& Ke^{-t/K\tau\log\tau} + \sum_{k=0}^{\lfloor t\rfloor\wedge
n} \p(S_{N+1} = n-k \& T_{N+2} > t + k\;\& N>0)\\
&+\p(T_1+T_2 > t)
\\
\le& Ke^{-t/K\tau\log\tau} + \sum_{k=0}^{\lfloor t\rfloor\wedge
n}\sum_{l=2}^{n-k} \p(S_l = n - k \; \&
T_{l+1} > t + k) + 2e^{-t/2\tau}\\
\le& Ke^{-t/K\tau\log\tau} + \sum_{k=0}^{\lfloor
t\rfloor}\sum_{l=2}^{n-k}\p(S_l = n - k)\p( T_{l+1}
> t + k)\\
\le& Ke^{-t/K\tau\log\tau} + 2(t+1)e^{-t/\tau} \le
Ke^{-t/K\tau\log\tau}.
\end{align*}

This implies that $\E(S_{N+2} - n)\le K\tau\log\tau \le K\tau^3/\E
T_2$, which proves (\ref{tofinish}). Thus (\ref{toprove}) is shown
and Theorem \ref{empirical_Markov} follows.
\end{proof}

\paragraph{Remark} Two natural questions to ask in regard to Theorem
\ref{empirical_Markov} is first whether the constant $K$ in front
of the expectation can be reduced to $1+ \eta$ (as in Massart's
Theorem \ref{Massart} or Theorem \ref{tail_estimate}) and second,
whether one can reduce the constant $K$ in the Gaussian part to
$2(1+\delta)$ (as in Theorem \ref{tail_estimate}).

\subsection{Another counterexample \label{section_counterexample}}

If we do not pay attention to constants, the main difference
between inequalities presented in the previous section and the
classical Bernstein's inequality for sums of i.i.d. bounded
variables is the presence of the additional factor $\log n$. We
would now like to argue that under the assumptions of Theorems
\ref{sums_Markov} and \ref{empirical_Markov}, this additional
factor is indispensable.

To be more precise, we will construct a Markov chain on a
countable state space, satisfying the assumptions of Theorem
\ref{sums_Markov} with $m = 1$ and such that for $\beta < 1$,
there is no constant $K$, such that
\begin{align}\label{hypothetical_inequality}
\p\Big(|f(X_1)+\ldots+f(X_n)| \ge t \Big) \le
K\exp\Big(-\frac{1}{K}\min\Big(\frac{t^2}{n
\Var(Z_1(f))},\frac{t}{\log^\beta n}\Big)\Big)
\end{align}
for all $n$ and all functions $f\colon \mathcal{S} \to \R$, with
$\|f\|_\infty \le 1$ and $\E_\pi f = 0$.

The state space of the chain will be the set
\begin{displaymath}
\mathcal{S} = \{0\}\cup\bigcup_{n=1}^\infty (\{n\}\times \{1,2,\ldots,n\}\times\{+1,-1\}).
\end{displaymath}
The transition probabilities are as follows
\begin{align*}
p_{(n,k,s),(n,k+1,s)} &= 1 \quad{\rm for}\quad n = 1,2,\ldots ,\; k = 1,2,\ldots,n-1,\; s = -1,+1\\
p_{(n,n,s),0} & = 1 \quad{\rm for}\quad n = 1,2,\ldots,\; s = -1,+1,\\
p_{0,(n,1,s)} & = \frac{1}{2A}e^{-n} \quad{\rm for}\quad n = 1,2,\ldots,\; s = -1,+1,\\
\end{align*}
where $A = \sum_{n=1}^\infty e^{-n}$. In other words, whenever a
''particle'' is at 0, it chooses one of countably many loops and
travels deterministically along it until the next return to $0$.
It is easy to check that this chain has a stationary distribution
$\pi$, given by
\begin{align*}
\pi_0 &= \frac{A}{A + \sum_{n=1}^\infty ne^{-n}},\\
\pi_{(n,i,s)} &= \frac{1}{2A}e^{-n}\pi_0.
\end{align*}

The chain satisfies the minorization condition
(\ref{minorization}) with $C = \{0\}$, $\nu(\{x\}) = p_{0,x}$,
$\delta = 1$ and $m = 1$. The random variable $T_1$ is now just
the time of the first visit to 0 and $T_2,T_3,\ldots$ indicate the
time between consecutive visits to 0. Moreover
\begin{displaymath}
\p(T_2 = n) = \frac{e^{-n}}{A},
\end{displaymath}
so $\|T_2\|_{\psi_1} < \infty$. If we start the chain from initial
distribution $\nu$, then $T_1$ has the same law as $T_2$, so $\tau
= \|T_2\|_{\psi_1} = \|T_1\|_{\psi_1}$.

Let us now assume that for some $\beta < 1$, there is a constant
$K$, such that (\ref{hypothetical_inequality}) holds. Since we
work with a fixed chain, in what follows we will use the letter
$K$ also to denote constants depending on our chain (the value of
$K$ may again differ at different occurrences).

We can in particular apply
(\ref{hypothetical_inequality}) to the function $f= f_r$ (where
$r$ is a large integer), given by the formula
\begin{displaymath}
f(0) = 0,\; f((n,i,s)) = s\ind{n \ge r}.
\end{displaymath}
We have $\E_{\pi} f_r = 0$. Moreover
\begin{displaymath}
\Var (Z_1(f_r)) = \sum_{n=r}^\infty n^2e^{-n}A^{-1} \le
Kr^2e^{-r}.
\end{displaymath}
Therefore (\ref{hypothetical_inequality}) gives
\begin{align}\label{hypothetical_inequality_2}
\p( |f_r(X_1) + \ldots + f_r(X_n)| \ge K(re^{-r/2}\sqrt{nt} +
t\log^\beta n)) \le e^{-t}
\end{align}
for $t \ge 1$ and $n\in \N$.

Recall that $S_i = T_1 + \ldots + T_i$. By Bernstein's inequality
(Lemma \ref{Bernstei_pis_1}), we have for large $n$,
\begin{align*}
\p(S_{\lceil n/(3\E T_2)\rceil} > n) &= \p(T_1 + \ldots + T_{\lceil n/(3\E T_2)\rceil} > n) \\
& = \p\Big(\sum_{i=1}^{\lceil n/(3\E T_2)\rceil} (T_i - \E T_i) > n - \lceil n/(3\E T_2)\rceil\E T_2\Big)\\
&\le \p\Big(\sum_{i=1}^{\lceil n/(3\E T_2)\rceil} (T_i - \E T_i) > n/2\Big)\\
&\le 2\exp\Big(-\frac{1}{K}\min
\Big(\frac{n^2}{n\|T_2\|_{\psi_1}^2},\frac{n}{\|T_2\|_{\psi_1}}\Big)\Big)
= 2e^{-n/K}.
\end{align*}
From the above estimate, for some integer $L$ and $n$ large enough, divisible by $L$,

\begin{align*}
\p\Big(&\Big|\sum_{i=0}^{n/L} Z_i(f_r)\Big| \ge 2K(re^{-r/2}\sqrt{nt} + t\log^\beta n)\Big) \\
\le&  2e^{-n/K} + \p\Big( \Big|\sum_{i=0}^{n/L} Z_i(f_r)\Big|\ge 2K(re^{-r/2}\sqrt{nt} + t\log^\beta n)\; \& S_{n/L+1} \le n \Big)\\
\le&   2e^{-n/K} + \p\Big(\Big|\sum_{i=0}^n f_r(X_i) \Big|\ge K(re^{-r/2}\sqrt{nt} + t\log^\beta n)\; \& S_{n/L+1} \le n \Big) \\
&+ \p\Big(\Big|\sum_{i=S_{n/L+1}+1}^n f_r(X_i)\Big|\ge K(re^{-r/2}\sqrt{nt} + t\log^\beta n)\; \& S_{n/L+1} \le n \Big)\\
\le& 2e^{-n/K} + e^{-t}\\
&+ \sum_{k\le n}\p\Big(\Big|\sum_{i=S_{n/L+1}+1}^n f_r(X_i)\Big|\ge K(re^{-r/2}\sqrt{nt} + t\log^\beta n)\; \& S_{n/L+1} = k\Big)\\
= &2e^{-n/K} + e^{-t}\\
&+\sum_{k\le n}\E\Big[\ind{S_{n/L+1} = k}\p\Big(\Big|\sum_{i=1}^{n-k} f_r(X_i)\Big|\ge K(re^{-r/2}\sqrt{nt} + t\log^\beta n)\Big)\Big]\\
\le&2e^{-n/K} + e^{-t} + e^{-t}\sum_{k\le n} \E\ind{S_{n/L+1}=k}
\le 2e^{-n/K} + 2e^{-t},
\end{align*}
where in the third and fourth inequality we used
(\ref{hypothetical_inequality_2}) and in the equality, the Markov
property.

For $n \simeq r^{-2}e^{r}$ and $t\ge 1$, we obtain
\begin{align} \label{hypothetical_estimate_3}
\p\Big(&|Z_0(f_r) + \ldots+ Z_{n/L}(f_r)| \ge K t\log^\beta n\Big)
\le 2e^{-t} + 2e^{-n/K}
\end{align}

On the other hand we have
\begin{align*}
\p(|Z_i(f_r)| \ge r) > \frac{1}{2A}e^{-r}.
\end{align*}
Therefore $\p(\max_{i\le n/L} |Z_i(f_r)| > r) \ge
2^{-1}\min(ne^{-r}/(2AL),1)$. Since $Z_i(f_r)$ are symmetric, by
Levy's inequality, we get
\begin{displaymath}
2\p\Big(|Z_0(f_r) + \ldots+ Z_{n/L}(f_r)| \ge r\Big) \ge
\frac{1}{2}\min(ne^{-r}/(2AL),1) \ge \frac{c}{r^2},
\end{displaymath}
whereas (\ref{hypothetical_estimate_3}) applied for $t =
K^{-1}r/\log^{\beta} n \ge K^{-1}r^{1-\beta} \ge 1$ gives
\begin{displaymath}
\p\Big(|Z_0(f_r) + \ldots+ Z_{n/L}(f_r)| \ge r\Big) \le 2e^{-
r^{1-\beta}/K} + 2e^{-e^r/(Kr^2)},
\end{displaymath}
which gives a contradiction.

\subsection{A bounded difference type inequality for symmetric
functions \label{boundeddifference}}

Now we will present an inequality for more general statistics of
the chain. Under the same assumptions on the chain as above (with
an additional restriction that $m=1$), we will prove a version of
the bounded difference inequality for symmetric functions (see
e.g. \cite{L} for the classical i.i.d. case).

Let us consider a measurable function $f \colon \mathcal{S}^n \to
\R$ which is invariant under permutations of arguments i.e.
\begin{align}\label{invariance}
f(x_1,\ldots,x_n) = f(x_{\sigma_1},\ldots,x_{\sigma_n})
\end{align}
for all permutations $\sigma$ of the set $\{1,\ldots,n\}$.

Let us also assume that $f$ is $L$-Lipschitz with respect to the
Hamming distance, i.e.
\begin{align}\label{Lipschitz_cond}
|f(x_1,\ldots,x_n) - f(y_1,\ldots,y_n)| \le L\#\{i\colon x_i\neq
y_i\}.
\end{align}

Then we have the following

\begin{theorem}\label{bounded_difference}
Let $X_1,X_2,\ldots$ be a Markov chain with values in
$\mathcal{S}$, satisfying the \textbf{Minorization condition}
with $m = 1$ and admitting a (unique) stationary distribution $\pi$. Assume also that $\|T_1\|_{\psi_1},\|T_2\|_{\psi_1} \le \tau$.
Then for every function $f \colon \mathcal{S}^n \to \R$,
satisfying (\ref{invariance}) and (\ref{Lipschitz_cond}), we have
\begin{displaymath}
\p(|f(X_1,\ldots,X_n) - \E f(X_1,\ldots,X_n)| \ge t) \le
2\exp\Big(-\frac{1}{K}\frac{t^2}{nL^2\tau^2}\Big)
\end{displaymath}
for all $t\ge 0$.
\end{theorem}

To prove the above theorem, we will need the following
\begin{lemma}
\label{lemma_symmet} Let $\varphi \colon \R \to \R$ be a convex
function and $G = f(Y_1,\ldots,Y_n)$, where $Y_1,\ldots,Y_n$ are
independent random variables with values in a measurable space
$\mathcal{E}$ and $f \colon \mathcal{E}^n\to \R$ is a measurable
function. Denote
\begin{displaymath}
G_i = f(Y_1,\ldots,Y_{i-1},\tilde{Y}_i,Y_{i+1},\ldots,Y_n),
\end{displaymath}
where $(\tilde{Y}_1,\ldots, \tilde{Y}_n)$ is an independent copy
of $(Y_1,\ldots,Y_n)$. Assume moreover that
\begin{displaymath}
|G - G_i| \le F_i(Y_i,\tilde{Y}_i)
\end{displaymath}
for some functions $F_i\colon \mathcal{E}^2 \to \R$, $i
=1,\ldots,n$. Then
\begin{equation}
\E \varphi(G - \E G) \le \E \varphi(\sum_{i=1}^n
\varepsilon_iF_i(Y_i,\tilde{Y}_i)), \label{symmet}
\end{equation}
where $\varepsilon_1,\ldots,\varepsilon_n$ is a sequence of
independent Rademacher variables, independent of $(Y_i)_{i=1}^n$
and $(\tilde{Y}_i)_{i=1}^n$.
\end{lemma}

\proof Induction with respect to $n$. For $n=0$ the statement is
obvious, since both the left-hand and the right-hand side of
(\ref{symmet}) equal $\varphi(0)$. Let us therefore assume that
the lemma is true for  $n-1$. Then, denoting by $\E_X$ integration
with respect to the variable $X$,
\begin{eqnarray*}
\E \varphi(G - \E G) & = & \E \varphi(G - \E_{\tilde{Y}_n}G_n + \E_{Y_n}G - \E G) \\
&\le& \E \varphi(G - G_n + \E_{Y_n}G - \E G) = \E \varphi(G_n - G+
\E_{Y_n}G - \E G)
\\ &=&
\E \varphi(\varepsilon_n|G-G_n| + \E_{Y_n}G - \E G) \\
&\le& \E \varphi(\varepsilon_nF_n(Y_n,\tilde{Y}_n) + \E_{Y_n}G -
\E G),
\end{eqnarray*}
where the equalities follow from the symmetry and the last
inequality from the contraction principle (or simply convexity of $\varphi$), applied conditionally
on $(Y_i)_i, (\tilde{Y_i})_i$. Now, denoting $Z = \E_{Y_n}G$, $Z_i
= \E_{Y_n}G_i$, we have for $i=1,\ldots,n-1$,
\begin{displaymath}
|Z - Z_i| = |\E_{Y_n}G - \E_{Y_n}G_i| \le \E_{Y_n}|G - G_i| \le
F_i(Y_i,\tilde{Y}_i),
\end{displaymath}
and thus for fixed $Y_n$,$\tilde{Y}_n$ and $\varepsilon_n$, we can
apply the induction assumption to the function $t \mapsto
\varphi(\varepsilon_nF(Y_n,\tilde{Y}_n) + t)$ instead of $\varphi$
and $\E_{Y_n}G$ instead of $G$, to obtain
\begin{displaymath}
\E \varphi(G - \E G) \le \E \varphi\left(\sum_{i=1}^n
F_i(Y_i,\tilde{Y}_i)\varepsilon_i\right).
\end{displaymath}
\qed

\begin{lemma}\label{psi_1_difference}
In the setting of Lemma \ref{lemma_symmet}, if for all $i$,
$\|F_i(Y_i,\tilde{Y_i})\|_{\psi_1} \le \tau$, then for all $t >
0$,
\begin{displaymath}
\p(|f(Y_1,\ldots,Y_n) - \E f(Y_1,\ldots,Y_n)| \ge t) \le
2\exp\Big(-\frac{1}{K}\min\Big(\frac{t^2}{n\tau^2},\frac{t}{\tau}\Big)\Big).
\end{displaymath}
\end{lemma}
\begin{proof}
For $p \ge 1$,
\begin{align*}
\|f(Y_1,\ldots,Y_n) - \E f(Y_1,\ldots,Y_n)\|_p &\le
\Big\|\sum_{i=1}^n\varepsilon_i F(Y_i,\tilde{Y}_i)\Big \|_p \le K(
\sqrt{p} \sqrt{n}\tau + p\tau),
\end{align*}
where the first inequality follows from Lemma \ref{lemma_symmet}
and the second one from Bernstein's inequality (Lemma
\ref{Bernstei_pis_1}) and integration by parts. Now, by the
Chebyshev inequality we get
\begin{displaymath}
\p(|f(Y_1,\ldots,Y_n) - \E f(Y_1,\ldots,Y_n)| \ge K(\sqrt{tn} +
t)\tau) \le e^{-t}
\end{displaymath}
for $t \ge 1$, which is up to the constant in the exponent
equivalent to the statement of the lemma (note that if we can
change the constant in the exponent, the choice of the constant in
front of the exponent is arbitrary, provided it is bigger than 1).
\end{proof}

\begin{proof}[Proof of Theorem \ref{bounded_difference}]
Consider a disjoint union
\begin{displaymath}
\mathcal{E} = \bigcup_{i=1}^\infty \mathcal{S}^i
\end{displaymath}
and a function $\tilde{f} \colon \mathcal{E}^n \to \R$ defined as
\begin{displaymath}
\tilde{f}(y_1,\ldots,y_n) = f(x_1,\ldots,x_n),
\end{displaymath}
where $x_i$'s are defined by the condition
\begin{align}\label{correspondence}
y_1 &= (x_1,\ldots,x_{t_1}) \in \mathcal{S}^{t_1}\nonumber\\
y_2 &= (x_{t_1+1},\ldots,x_{t_1 + t_2}) \in \mathcal{S}^{t_2}\nonumber\\
\ldots\nonumber\\
y_n & = (x_{t_1+\ldots+t_{n-1}+1},\ldots,x_{t_1+\ldots+t_n}) \in
\mathcal{S}^{t_n}.
\end{align}
Let now $T_1,\ldots,T_n$ be the regeneration times of the chain
and set
\begin{displaymath}
Y_i = (X_{T_1+\ldots+T_{i-1}+1},\ldots,X_{T_1+\ldots +
T_i})
\end{displaymath}
for $i = 1,\ldots,n$ (we change the enumeration with
respect to previous sections, but there is no longer need to
distinguish the initial block). Then $Y_1,\ldots,Y_n$ are
independent $\mathcal{E}$-valued random variables (recall the
assumption $m=1$). Moreover we have
\begin{displaymath}
f(X_1,\ldots,X_n) = \tilde{f}(Y_1,\ldots,Y_n).
\end{displaymath}
Let now $\tilde{Y_1},\ldots,\tilde{Y}_n$ be an independent copy of
the sequence $Y_1,\ldots,Y_n$. Define $G$ and $G_i$ like in Lemma
\ref{lemma_symmet} (for the function $\tilde{f}$). Define also
$\tilde{T}_i = j$ iff $\tilde{Y}_i \in \mathcal{S}^j$ and let
$\tilde{X}_{i,1},\ldots,\tilde{X}_{i,T_1+\ldots+T_{i-1}+\tilde{T}_i+T_{i+1}+\ldots
+T_n}$ correspond to
$Y_1,\ldots,Y_{i-1},\tilde{Y}_i,Y_{i+1},\ldots,Y_n$ in the same
way as in (\ref{correspondence}). Let us notice that we can
rearrange the sequence $(\tilde{X}_{i,1},\ldots,\tilde{X}_{i,n})$
in such a way that the Hamming distance of the new sequence from
$(X_1,\ldots,X_n)$ will not exceed $\max(T_i,\tilde{T}_i)$. Since
the function $f$ is invariant under permutation of arguments and
$L$-Lipschitz with respect to the Hamming distance, we have
\begin{displaymath}
|G - G_i| \le L \max(T_i,\tilde{T}_i) =: F(Y_i,\tilde{Y}_i).
\end{displaymath}
Moreover, $\|F(Y_i,\tilde{Y}_i)\|_{\psi_1} \le 2L\tau$, so by
Lemma \ref{psi_1_difference}, we obtain
\begin{align*}
&\p(|f(X_1,\ldots,X_n) - \E f(X_1,\ldots,X_n)| \ge t) \\
&=\p(|\tilde{f}(Y_1,\ldots,Y_n) - \E \tilde{f}(Y_1,\ldots,Y_n)|
\ge t) \le
2\exp\Big(-\frac{1}{K}\min\Big(\frac{t^2}{nL^2\tau^2},\frac{t}{L
\tau}\Big)\Big).
\end{align*}
But from Jensen's inequality and (\ref{Lipschitz_cond}) it follows
that $|f(X_1,\ldots,X_n) - \E f(X_1,\ldots,X_n)| \le L n$, thus
for $t > L n$, the left hand side of the above inequality is equal
to $0$, whereas for $t \le L n$, the inequality $\tau > 1$ gives
\begin{displaymath}
\frac{t^2}{nL^2\tau^2} \le \frac{t}{L\tau},
\end{displaymath}
which proves the theorem.
\end{proof}

\subsection{A few words on connections with other results
\label{discussion}} First we would like to comment on the
assumptions of our main theorems, concerning Markov chains. We
assume that the Orlicz norms $\|T_1\|_{\psi_1}$ and
$\|T_2\|_{\psi_1}$ are finite, which is equivalent to existence of
a number $\kappa > 1$, such that
\begin{displaymath}
\E_\xi \kappa^{T_1} < \infty, \; \E_\nu \kappa^{T_1} < \infty,
\end{displaymath}
where $\xi$ is the initial distribution of the chain and $\nu$ --
the minorizing measure from condition (\ref{minorization}). This
is true for instance if $m=1$ and the chain satisfies the
drift condition, i.e. if there is a measurable function $V \colon
\mathcal{S} \to [1,\infty)$, together with constants $\lambda< 1$
and $K < \infty$, such that
\begin{align*}
PV(x) = \int_\mathcal{S} V(y)P(x,dy) \le \left\{
\begin{array}{rcl}
\lambda V(x) &\;{\rm for} \; & x \notin C,\\
K&\;{\rm for} \; & x \in C
\end{array}\right.
\end{align*}
and $V$ is $\xi$ and $\nu$ integrable (see e.g. \cite{B},
Propositions 4.1 and 4.4, see also \cite{RR}, \cite{MT}). For $m >
1$ one can similarly consider the kernel $P^m$ instead of $P$
(however in this case our inequalities are restricted to averages
of real valued functions as in Theorem \ref{sums_Markov}). Such
drift conditions have gained considerable attention in the Markov
Chain Monte Carlo theory as they imply geometric ergodicity of the
chain.

Concentration of measure inequalities for general functions of
Markov chains were investigated by Marton \cite{M}, Samson
\cite{S} and more recently by Kontorovich and Ramanan \cite{K}.
They actually consider more general mixing processes and give
estimates on the deviation of a random variable from the mean or
median in terms of mixing coefficients. When specialized to Markov
chains, their estimates yield inequalities in the spirit of
Theorem \ref{bounded_difference} for general (non-necessarily
symmetric) functions of uniformly ergodic Markov chains (see
\cite{MT}, Chapter 16 for the definition). To obtain their
results, Marton and Samson used transportation inequalities,
whereas Kontorovich's and Ramanan's approach was based on
martingales. In all cases the bounds include sums of expressions
of the form
\begin{displaymath}
\sup_{x,y\in\mathcal{S}} \|P^i(x,\cdot) - P^i(y,\cdot)\|_{\rm TV},
\end{displaymath}
where $P^i$ is the $i$ step transition function of the chain.
These results are not well suited for Markov chains which are not
uniformly ergodic (like the chain in Section
\ref{section_counterexample}), since for such chains the summands
are bounded from below by a constant (which spoils the dependence
on $n$ in the estimates). It would be interesting to know if in
results of this type, the supremum of the total variation
distances can be replaced by some other norm, for instance a kind
of average. This would allow to extend the estimates to some
classes of non-uniformly ergodic Markov chains.

\paragraph{}
Inequalities of the bounded difference type for sums
$f(X_1)+\ldots+ f(X_n)$ where $X_i$'s form a uniformly ergodic
Markov chain were also obtained by Glyn and Ormoneit \cite{G}.
Their method was to analyze the Poisson equation associated with
the chain. Their result has been complemented by an information
theoretic approach in Kontoyiannis et al. \cite{Ky}.

\paragraph{}
Estimates for sums, in terms of variance, appeared in the work by
Samson \cite{S}, who presents a result for empirical processes of
uniformly ergodic chains. He gives a real concentration inequality
around the mean (and not just a tail bound as in Theorem
\ref{empirical_Markov}). The coefficient responsible for the
subgaussian behavior of the tail is $\E\sum_{i=1}^n \sup_f
f(X_i)^2$. Replacing it with $V = \E \sup_f \sum_i f(X_i)^2$
(which would correspond to the original Talagrand's inequality) is
stated in Samson's work as an open problem, which to our best
knowledge has not been yet solved. Additionally, in Samson's
estimate there is no $\log n$ factor, which is present in Theorems
\ref{sums_Markov} and \ref{empirical_Markov}. Since we have shown
that in our setting this factor is indispensable, we would like to
comment on the differences between the results by Samson and ours.

Obviously, the first difference is the setting. Although
non-uniformly ergodic chains satisfy our assumptions
$\|T_1\|_{\psi_1},\|T_2\|_{\psi_1} < \infty$, the Minorization
condition may not hold for them with $m = 1$, which restricts our
results to linear statistics of the chain (Theorem
\ref{sums_Markov}). However, there are many examples of
non-uniformly ergodic chains, for which one cannot apply Samson's
result but which satisfy our assumptions. Such chains have been
considered in the MCMC theory.

When specialized to sums of real variables, Samson's result can be
considered a counterpart of the Bernstein inequality, valid for
uniformly ergodic Markov chains. The subgaussian part of the
estimate is controlled by $\sum_{i=1}^n \E f(X_i)^2$, which can be
much bigger than the asymptotic variance and therefore does not
reflect the limiting behaviour of $f(X_1)+\ldots + f(X_n)$.
Consider for instance a chain consisting of the origin connected
with finitely many loops in which, similarly as in the example
from Section \ref{section_counterexample}, the randomness appears
only at the origin (i.e. after the choice of the loop the particle
travels along it deterministically until the next return to the
origin). Then, one can easily construct a function $f$ with values
in $\{\pm 1\}$, centered with respect to the stationary
distribution and such that its asymptotic variance is equal to
zero, whereas $\sum_{i=1}^n \E f(X_i)^2 = n$ for all $n$ (it
happens for instance if the sum of the values of $f$ along each
loop vanishes). In consequence, $n^{-1/2}(f(X_1)+\ldots+f(X_n))$
converges weakly to the Dirac mass at $0$ and we have
\begin{displaymath}
\p(|f(X_1)+\ldots+f(X_n)| \ge \sqrt{n} t) \to 0
\end{displaymath}
for all $t > 0$, which is not recovered by Samson's estimate. One
can also construct other examples of similar flavour, in which the
asymptotic variance is nonzero but is still much smaller than $\E
\sum_{i=1}^n f(X_i)^2$.

On the other hand Samson's results do not require the condition
$\E_\pi f = 0$ and (as already mentioned) in the case of empirical
processes they provide a two sided concentration around the mean.

As for the $\log n$ factor, at present we do not know if at the
cost of replacing the asymptotic variance with $\sum_{i=1}^n \E
f(X_i)^2$ one can eliminate it in our setting.

\paragraph{}
Summarizing, our inequalities, when compared to known results have
both advantages and disadvantages. On the one hand, when
specialized to uniformly ergodic Markov chains, they do not
recover the full generality or strength of previous estimates (for
instance Theorem \ref{bounded_difference} is restricted to
symmetric statistics and $m=1$), on the other hand they may be
applied to Markov chains arising in statistical applications,
which are not uniformly ergodic (and therefore beyond the scope of
the estimates presented above). Another property, which in our
opinion, makes the estimates of Theorems \ref{sums_Markov} and
\ref{empirical_Markov} interesting (at least from the theoretical
point of view) is the fact that for $m=1$, the coefficient
responsible for the Gaussian level of concentration corresponds to
the variance of the limiting Gaussian distribution.

\paragraph{Acknowledgements} The author would like to thank Witold
Bednorz and Krzysztof {\L}atuszy\'{n}ski for their useful comments
concerning the results presented in this article as well as the
anonymous Referee, whose remarks helped improve their
presentation.

\end{document}